\documentclass[11pt,draft]{amsart}
\usepackage{dsfont}
\usepackage{amsfonts}
\usepackage{mathrsfs}
\usepackage{amssymb}
\usepackage{amsmath}
\usepackage{amsfonts}
\usepackage{amsfonts,enumerate}
\usepackage{pifont}
\usepackage{enumerate}
\usepackage{graphics}
\usepackage{verbatim}
\usepackage{amsthm}
\usepackage{latexsym,bm}
\usepackage{color}
\numberwithin{equation}{section}
\newtheorem{theorem}{Theorem}[section]
\newtheorem{corollary}[theorem]{Corollary}
\newtheorem{lemma}[theorem]{Lemma}
\newtheorem{proposition}[theorem]{Proposition}

\theoremstyle{definition}
\newtheorem{remark}[theorem]{Remark}

\newcommand{\B}{{\mathcal B}}

\begin{document}

\title[]{Vector valued $q$-variation for differential operators and semigroups II}

\author{Guixiang Hong*}
\address{Instituto de Ciencias Matem\'aticas,
CSIC-UAM-UC3M-UCM, Consejo Superior de Investigaciones
Cient\'ificas, C/Nicol\'as Cabrera 13-15. \newline 28049, Madrid. Spain.\\
\emph{E-mail address: guixiang.hong@icmat.es}}
\thanks{*\ Corresponding author}

\author{Tao Ma}
\address{School of Mathematics and Statistics, Wuhan University, Wuhan 430072, China \\ \emph{E-mail address: tma.math@whu.edu.cn}}

\thanks{\small {{\it MR(2010) Subject Classification}.} Primary
47B38, 46E40; Secondary 47A35, 42B25.}

\maketitle

\begin{abstract}
In this paper, we establish UMD lattice-valued variational inequalities for differential operators, ergodic averages and analytic semigroups. These results generalize, on the one hand some scalar-valued variational inequalities in ergodic theory, on the other hand Xu's very recent result on UMD lattice-valued maximal inequality. As a consequence, we deduce the jump estimates and  obtain quantitative information on the rate of the pointwise convergence.
\end{abstract}

\section{Introduction}
Let $(\Omega,\mu)$ be a measure space and $\B$ be a Banach space. A submarkovian $C_0$ semigroup $(T_t)_{t\geq0}$ acting on $L^p(\Omega)$ extends to a semigroup of operators on $L^p(\Omega;\B)$. Cowling and Leinert showed in \cite{CoLe11} that
\begin{align}\label{motivation1}
\|T_tf(\omega)-f(\omega)\|\rightarrow 0,\;\mathrm{a.e.}\;\omega\in\Omega,\;\mathrm{as}\;t\rightarrow0^+
\end{align}
for any Banach space $\B$ and any $f\in L^p(\Omega;\B)$ with $1<p<\infty$, where $\|\cdot\|$ means taking $\B$-norm. The predecessor of this result is the one \cite{Tag09} by Taggart where the convergence was shown only for Banach spaces having UMD property.

When $\B$ is a Banach lattice of measurable functions on a measure space $(\Sigma,\nu)$,  $\B$-valued functions on $\Omega$ can be viewed as scalar-valued functions on $\Omega\times\Sigma$. Different from previous paper \cite{HoMa}, in this paper, we are concerned with the question whether we have
\begin{align}\label{motivation2}
T_tf(\omega,\sigma)\rightarrow f(\omega,\sigma),\;\mathrm{a.e.}\;(\omega,\sigma)\in\Omega\times\Sigma,\;\mathrm{as}\;t\rightarrow0^+
\end{align}
for any Banach lattice $\B$ and any $f\in L^p(\Omega;\B)$ with $1<p<\infty$?
A priori the answer would seem to be no, since then Stein's principle \cite{Ste61} would imply a vector-valued maximal inequality would hold for all Banach lattices, which is false since then all Banach lattices have H.L. property which contradict with the fact that $\ell^1$ does not have H.L. property (see \cite{MaTo}).

On the other hand, Xu \cite{Xu} has recently established a Banach lattice-valued maximal inequality which implies
the pointwise convergence (\ref{motivation2}) via Banach principle under the condition that $\B$ is an UMD lattice.

In this paper, we prove that if $\B$ is an UMD lattice, some Banach lattice-valued variational inequalities still hold which immediately imply the pointwise convergence of the underlying family of operators without using the Banach principle. Moreover, the variational inequalities provide us the quantitative information of the rate of the convergence.

The scalar-valued variational inequalities have been studied a lot in probability, ergodic theory and harmonic analysis. The first variational inequality was proved by L\'epingle \cite{Lep76} for martingales which improves the classical Doob maximal inequality (see also \cite{PiXu88} for a different approach and related results). Thirteen years later, Bourgain in \cite{Bou89} proved the variational inequality for the ergodic averages of a dynamic system, which has inaugurated a new research direction in ergodic theory and harmonic analysis.  Bourgain's work was considerably improved and extended to many other operators in ergodic theory. For instance, Jone et al in \cite{JKRW98} and \cite{JRW03}  established variational inequalities for differential operators and ergodic averages of measure-preserving invertible transforms. On the other hand, almost in the same period, variational inequalities have been studied in harmonic analysis too. The first works on this subject are \cite{CJRW00} and \cite{CJRW03} where Campbell et al proved the variational inequalities for singular integrals. Since then variational inequalities for different kinds of operators in harmonic analysis have been built; see e.g. \cite{DMT12}, \cite{JoWa04}, \cite{Mas}, \cite{MaTo12},  \cite{MaTo}, \cite{OSTTW12}.

Among the works, we would like to mention the results obtained recently by (Ma/Torrea/Xu \cite{MTX1}) Le Merdy/Xu \cite{LeXu2}, where (weighted) variational inequalities for (differential operator) semigroups were established, as well as Xu \cite{Xu} where vector-valued $H^\infty$ functional calculus for analytic semigorups on vector-valued $L^p$ spaces was built. The results or idea in these papers will be exploited in the present paper.

To state our results, we need to recall the definition of the $q$-variation. Give a sequence $(a_n)_{n\geq0}$ of complex number and a number $1\leq q<\infty$, the  $q$-variation norm is defined as
$$\|(a_n)_{n\geq0}\|_{v_q}=\sup\{(|a_{n_0}|^q+\sum_{k\geq1}|a_{n_{k}}-a_{n_{k-1}}|^{q})^{\frac{1}{q}}\}$$
where the supremum runs over all increasing sequences $(n_k)_{k\geq0}$ of integers. It is clear that the set of $v_q$ of all sequences with a finite  $q$-variation is a Banach space for the norm $v_q$.  A continuous analog of $v_q$ is defined as follows. Given a family $(a_t)_{t>0}$ of complex numbers, define
$$\|(a_t)_{t>0}\|_{V_q}=\sup\{(|a_{t_0}|^q+\sum_{k\geq1}|a_{t_{k}}-a_{t_{k-1}}|^{q})^{\frac{1}{q}}\}$$
where the supremum runs over all increasing sequences $(t_k)_{k\geq0}$ of positive real numbers. Then we define $V_q$ to be the Banach space of all $(a_t)_{t>0}$ with $V_q$-norm finite.

The main result of this paper is stated as follows.  Let $(T_t)_{t\geq0}$ be a bounded analytic semigroup on $L^p(\Omega)$ with $1<p<\infty$ and assume that $T_t$ is a contractively regular for any $t\geq0$. Let $2<q<\infty$ and $\B$ be an UMD lattice having Fatou property. Then for any $f\in L^p(\Omega;\B)$,  the family $(T_tf)_{t>0}$ belongs to $V_q$ for a.e. $(\omega,\lambda)\in\Omega\times\Sigma$ and we have an estimate
\begin{align}\label{lattice analytic semigroup continuous0}
\left\|\|(\omega,\lambda)\rightarrow (T_tf(\omega,\lambda))_{t>0}\|_{V_q}\right\|_{L^p(\Omega;\B)}\leq C_{p,q} \|\|f\|\|_p,\;\forall f\in L^p(\Omega;\B).
\end{align}

Some comments on this result are in order. First of all, the semigroup under consideration is much general than the one in Cowling and Leinert, since analytic and contractively regular is a much weaker condition than submarkovian (see e.g. \cite{Gol85} and \cite{Paz83}). Secondly, inequality (\ref{lattice analytic semigroup continuous0}) allows us to conclude the convergence (\ref{motivation2}) and the speed of the convergence, see Section 5 for the statements; Thirdly, this result can be regarded as a Banach lattice-valued version of Corollary 4.5 in \cite{LeXu2} by Le Merdy and Xu; Finally, this result immediately implies the main result-Theorem 2 of \cite{Xu} by Xu in the radial case.

A priori, the kind of Banach lattice valued inequality as (\ref{lattice analytic semigroup continuous0}) is difficult to deal with, since it is well known that the Banach lattice valued inequalities are closely related to the weighted norm inequality in the harmonic analysis, and there is lack of the theory of weighted norm in the general ergodic setting. However, the arguments used in the scalar valued case \cite{LeXu2} is very powerful, where pointwise estimates are used, so that the general pattern can be adapted to the vector-valued case with more efforts. More precisely, as in \cite{LeXu2} we will deduce inequality (\ref{lattice analytic semigroup continuous0}) from a similar estimate for discrete semigroups using an approxiamtion argument based on the semigroup property, which in turn relys on a similar estimate for ergodic averages $M(T_t)$'s (see Section 3 for the definition)
\begin{align}\label{lattice ergodic average0}
\|(M_n(T_t)f)_{n\geq0}\|_{L^p(\B(v_q))}\|\leq C_{p,q}\|\|f\|\|_{p},\;\forall f\in L^p(\Omega;\B)
\end{align}
and  Banach lattice-valued Littlewood-Paley inequality
\begin{align}\label{lattice littlewood-paley0}
\left\|\left(\sum^{\infty}_{n=0}\frac{1}{n+1}|(n+1)T_t^n(T_t-I)(f)|^{2}\right)^{\frac{1}{2}}\right\|_{L^p(\B)}\leq C_p\|\|f\|\|_p,\;\forall f\in L^p(\Omega;\B).
\end{align}

The estimate (\ref{lattice ergodic average0}) will be proved through transference technique from the particular case, that is, the ergodic averages of the shift on $\mathbb{Z}$,
\begin{align}\label{lattice differential operator on Z10}
\|(A_nf)_{n\geq0}\|_{L^p(\B(v_q))}\|\leq C_{p,q}\|\|f\|\|_{p},\;\forall f\in L^p(\mathbb{Z}; \B),
\end{align}
for which we exploit a ergodic variant of Rubio's extrapolation theorem \cite{Rub86}, based on the weighted norm inequality established in \cite{MTX1}. This approach is very different from the one used for the vector-valued variational inequality considered in \cite{HoMa}, where all the arguments in the scalar-valued case should be adapted.

The estimate (\ref{lattice littlewood-paley0}) will be shown through approximation arguments based on the semigroup property and a continuous version of (\ref{lattice littlewood-paley0}) established by Xu in \cite{Xu} using vector-valued $H^\infty$ functional calculus. We also explain shortly how one can adapt the arguments in \cite{LeXu1} to deduce inequality  (\ref{lattice littlewood-paley0}) for any single analytic operator (replacing $T_t$ by $T$) from  Xu's vector-valued $H^\infty$ functional calculus. 

Our paper is organized as follows. In Section 2, we  prove inequality (\ref{lattice differential operator on Z10}) and the related results. In section 3, by transference principle, we show variational inequality (\ref{lattice ergodic average0}) for ergodic averages. In section 4, the variational inequality for the semigroups itself will be proved, based on (\ref{lattice ergodic average0}) and Littlewood-Paley inequality (\ref{lattice littlewood-paley0}). In the last section, we establish individual (pointwise) ergodic theorems and the quantitative formulation of convergence.

\section{Differential operators}
Let $\B$ be a Banach lattice of measurable functions on a measure space $(\Sigma,\nu)$. For $\B$-valued function $f$ defined on $\mathbb{Z}$, we define the differential operators for any integer $n\geq0$ as
$$A_nf(j)=\frac{1}{n+1}\sum^{n}_{k=0}f(j+k).$$

The main result of this section is stated as follows.
\begin{theorem}\label{thm: lattice differential operator on Z1}
Let $\B$ be an UMD lattice. Let $1<p<\infty$ and  $2<q<\infty$. Then there exists an absolute constant $C_{p,q}$ such that
\begin{align}\label{lattice differential operator on Zd}
\|(A_nf)_{n\geq0}\|_{L^p(\B(v_q))}\leq C_{p,q}\|f\|_{L^p(\B)}, \;\forall f\in L^p(\mathbb{Z}; \B).
\end{align}
\end{theorem}

\begin{remark}
When $\B=\ell^r$ with $1<r<\infty$, Theorem \ref{thm: lattice differential operator on Z1} is the unweighted case of Theorem 4.1 in \cite{MTX1} which is shown by the extrapolation theorem.
\end{remark}

In the case $\mathbb{R}^n$, a deep result by Rubio-Theorem 5 of \cite{Rub86} says that UMD lattice-valued inequality can be deduced from the weighted norm inequality. We show below a similar result still holds in the case $\mathbb{Z}$. Then Theorem \ref{thm: lattice differential operator on Z1} will follow immediately from this result since the weighted norm inequality has been built in \cite{MTX1}.

\begin{proposition}\label{pro: rubio characterization Z1}
Let $1<p<\infty$, $\B$ be an UMD lattice, and let $S$ be a not necessarily linear operator which is bounded in $L^p(\mathbb{Z}, w)$ for all $w\in A_p$. Then $S(f)(n,\lambda)=S(f(\cdot,\lambda))(n)$ is bounded in $L^p(\mathbb{Z}; \B)$.
\end{proposition}

As argued by Rubio, we need a lemma (Lemma 1 in \cite{Rub86}), which finds a deep connection between vector-valued extension and weighted norm estimate. We state it here for convenience.

\begin{lemma}\label{rubio lemma}
Let $Y$ be a superreflexive Banach lattice of measurable functions in $(\Omega,dw)$, and let $T: Y\rightarrow Y$ be a bounded linear operator. If $Y$ is a $p$-convex, the following two statements are equivalent:
\begin{enumerate}[\rm (i)]
\item $T$ has an extension $\tilde{T}: Y(\ell^p)\rightarrow Y(\ell^p)$, i.e.
$$\left\|\left(\sum_j|Ty_j|^p\right)^{\frac{1}{p}}\right\|_Y\leq C_1\left\|\left(\sum_j|y_j|^p\right)^{\frac{1}{p}}\right\|_Y;$$
\item For every $u\in(Y^p)^*$, there exists $w\in(Y^p)^*$ such that
$$u\leq w,\; \|w\|_{(Y^p)^*}\leq2\|u\|_{(Y^p)^*}$$ and
$$\left\|\int_{\Omega}|Ty|^pu\right\|^{\frac{1}{p}}\leq C_2\left\|\int_{\Omega}|y|^pw\right\|^{\frac{1}{p}}$$
where $Y^p$ is the $p$-convexification of $Y$.
\end{enumerate}
\end{lemma}

\begin{remark}
The implication from $\mathrm{(ii)}$ to $\mathrm{(i)}$ consists of an elementary application of the duality between $Y^p$ and $(Y^p)^*$, hence the operator $T$ is not necessarily  linear.
\end{remark}

\begin{proof}
Since UMD implies superreflexivity, every UMD lattice $\B$ is $p_0$-convex for certain $1<p_0<\infty$. Fix $1<p<p_0$. Since $\B(\ell^p)$ is still an UMD Banach space, by Theorem 2.8 of \cite{BeGi87}, the discrete Hilbert transform, defined as
$$H_{disc}((a_n)_{n})(m)=\sum_{n\neq m}\frac{a_n}{n-m}$$
is bounded in $L^{rp}(\mathbb{Z}; \B(\ell^p))$ for all $1<r<\infty$. Then apply Lemma \ref{rubio lemma} to $Y=L^{rp}(\mathbb{Z}; \B)$ which is obviously $p$-convex, we have
$$\sum_{n\in\mathbb{Z}}\int_{\Sigma}|H_{disc}f(n,\lambda)|^pw(n,\lambda)d\nu(\lambda)\leq C\sum_{n\in\mathbb{Z}}\int_{\Sigma}|f(n,\lambda)|^pw(n,\lambda)d\nu(\lambda)$$
for any $w\in (Y^p)^*$ and any $f\in L^p(\mathbb{Z};\B)$. In particular, apply this estimate for the function $f(n,\lambda)=a(n)\chi_{E}(\lambda)$, where $a$ is a finite sequence and $E$ is a subset of $\Sigma$ of finite measure. Then $H_{disc}f(n,\lambda)=\chi_{E}(\lambda)H_{disc}(a)(n)$, and we obtain
$$\sum_{n\in\mathbb{Z}}|H_{disc}(a)(n)|^pw(n,\lambda)\leq C\sum_{n\in\mathbb{Z}}|a(n)|^pw(n,\lambda) ,\;\mathrm{a.e.}\;\lambda\in\Sigma,$$
which implies, by Theorem 10 in \cite{HMW73}, $w(\cdot,\lambda)$ belong to $A_p$ for a.e. $\lambda\in\Sigma$. Then by our assumption, $S$ is bounded in  $L^p(\mathbb{Z}, w)$ for all $w\in A_p$, we deduce that for any $w\in (Y^p)^*$,
$$\sum_{n\in\mathbb{Z}}\int_{\Sigma}|Sf(n,\lambda)|^pw(n,\lambda)d\nu(\lambda)\leq C\sum_{n\in\mathbb{Z}}\int_{\Sigma}|f(n,\lambda)|^pw(n,\lambda)d\nu(\lambda).$$
Use again Lemma \ref{rubio lemma} (noting the remark afterwards), we obtain
$$\|Sf\|_{L^{rp}(\B)}\leq C\|f\|_{L^{rp}(\B)},$$
which implies the desired result by noting that $p$ varies between $1<p<p_0$ and $r$ varies between $1<r<\infty$.
\end{proof}

From the proof, a high dimensional version of Proposition \ref{pro: rubio characterization Z1}  would also be true if we could find an operator $T$ on $\mathbb{Z}^d$ playing the role of $H_{disc}$ on $\mathbb{Z}$ such that
\begin{enumerate}[(i)]
\item For any UMD lattice $\B$ and any $1<p<\infty$, we have $T$ is bounded on $L^p(\mathbb{Z}^d;\B)$;
\item If $T$ is bounded on $L^p(\mathbb{Z}^d,w)$ for all $p>1$ but close to $1$ then $w$ must belong to the  $A_p$-weights class on $\mathbb{Z}^d$.
\end{enumerate}

These results will appear elsewhere. On the other hand, in the case $\mathbb{R}^d$, the operator $\sum^d_{j=1}R_j$ with $R_j$'s Riesz transforms satisfies the condition (i) (see e.g. \cite{Zim89}) and (ii) (see e.g. P. 321 in \cite{Gra}), hence similar UMD lattice-valued variational inequalities as  inequality \ref{lattice differential operator on Zd} are true since the weighted norm inequalities for $q$-variation on $\mathbb{R}^d$ has been established recently by Ma, Torrea and Xu  \cite{MTX2}.

At the end of this section, we would like to show that the result in the case $\mathbb{R}$ can be deduced from the one in the case $\mathbb{Z}$ using an approximation argument without using the weighted theory and Rubio's result. For any positive real number $t>0$, we define
$$A_tf(s)=\frac{1}{t}\int^t_{0}f(r+s)dr,\forall f\in L_{1,loc}(\mathbb{R};\B).$$

\begin{corollary}\label{cor: lattice differential operator on R1}
Let $2<q<\infty$, $\B$ be an UMD lattice having Fatou property  and let $f\in L^p(\Omega;\B)$. Then for a.e. $(\omega,\lambda)\in(\Omega\times\Sigma)$, the family $((A_tf)(\omega,\lambda))_{t\geq0}$ belongs to $\B(V_q)$ and
$$\left\|(\omega,\lambda)\rightarrow\|((A_tf)(\omega,\lambda))_{t\geq0}\|_{\B(V_q)}\right\|_p\leq C_{p,q}\|\|f\|\|_p.$$
\end{corollary}

To prove this Corollary, we need the following lemma.

\begin{lemma}\label{lem: approximation properties 4}
Let $(f_t)_{t>0}$ be a family of $L^p(\Omega; \B)$ with $\B$ having Fatou property and assume that:
\begin{enumerate}[\rm(i)]
\item For a.e. $(\omega,\lambda)\in \Omega\times\Sigma$, the function $t\rightarrow f_t(\omega,\lambda)$ is continuous on $(0,\infty)$;
\item There exists a constant $C>0$ such that whenever $t_0<t_1<\dotsm<t_N$ is a finite increasing sequence of positive real numbers, we have
$$\|(f_{t_0},f_{t_1},\dotsm,f_{t_N})\|_{L^p(\Omega; \B(v_q))}\leq C.$$
\end{enumerate}
Then $(f_t(\omega,\lambda))_{t>0}$ belongs to $V_q$ for a.e. $(\omega,\lambda)\in\Omega\times\Sigma$, the function $(\omega,\lambda)\rightarrow \|(f_t(\omega,\lambda))_{t>0}\|_{V_q}$ belongs to $L^p(\Omega;\B)$ and
$$\left\|(\omega,\lambda)\rightarrow \|(f_t(\omega,\lambda))_{t>0}\|_{V_q}\|\right\|_{L^p(\B)}\leq C.$$
\end{lemma}
This lemma can be proved using the same argument for Lemma 2.2 in \cite{LeXu2}, since $F$ has Fatou property. We refer to \cite{LiTz79} for more information on Banach lattices and related properties.

Now we are at a position to give the proof of Corollary \ref{cor: lattice differential operator on R1}.

\begin{proof}
First of all, note that the function $t\rightarrow (A_tf)(\omega,\lambda) $ is continuous for a.e. $(\omega,\lambda)\in(\Omega\times\Sigma)$. Use the strong continuity of $A_t$ in $L^p(\mathbb{R};\B)$ (see e.g. \cite{Tag09}), for fixed positive numbers $t_0<t_1<\dotsm<t_N$ and $\varepsilon$, there exist $\eta>0$ and positive integers $n_0,\dotsm, n_N$ such that
$$\|A_{t_k}(f)-A_{n_k}(f(\eta\cdot))\|_{L^p(\B)}<\varepsilon\;\forall 0\leq k\leq N.$$
Hence by  Theorem  \ref{lattice ergodic average} and a limit argument, we deduce that
$$\|A_{t_0}f, A_{t_1}f,\dotsm, A_{t_m}f\|_{L^p(\B(v_q))}\leq C_{p,q}\|\|f\|\|_p.$$
The desired result therefore follows from Lemma \ref{lem: approximation properties 4}.
\end{proof}

\section{Ergodic averages}
It is well known that many results in harmonic analysis on $\mathbb{Z}$ can be transferred to the ergodic theory. In this section, we will show that the vector-valued $q$-variation for the ergodic averages is also bounded on $L^p(\B)$ with $1<p<\infty$.  Let $(\Omega_1,\mu)$ and $(\Omega_2,\mu)$ be two measure spaces.  An operator $T: L^p(\Omega_1)\rightarrow L^p(\Omega_2)$ is called contractively regular if the following inequality holds
$$\left\|\sup_{k\geq1}|T(x_k)|\right\|_{L^p(\Omega_1)}\leq \left\|\sup_{k\geq1}|x_k|\right\|_{L^p(\Omega_2)}.$$
for any finite sequence $(x_k)_{k\geq1}$ in $L^p(\Omega_1)$. Any contractively regular operator $T$ can be extended to a contractive operator on the Bochner space $L^p(\Omega_1; E)$ for any Banach space $E$, i.e.
\begin{align}\label{extension of T}
\|T\otimes I_{E}: L^p(\Omega_1;E)\rightarrow L^p(\Omega_2;E)\|\leq1.
\end{align}
In this paper, any extension of any operator $S$ will be still denoted by $S$ when no confusion occurs.

Obviously, positive contractions are regular. On the other hand, it is easy to check that if $T$ is a contraction on $L^1(\Omega_1)$ and $L^{\infty}(\Omega_1)$, then $T$ is contractively regular on  $L^p(\Omega_1)$. We refer to \cite{Mey91}, \cite{Pel76} and \cite{Pis94} for more details and complements.

 Let $T: L^p(\Omega)\rightarrow L^p(\Omega)$ be a contractively regular operator. For any integer $n\geq0$, the ergodic averages of $T$ is defined as
$$M_n(T)=\frac{1}{n+1}\sum^{n}_{k=0}T^k.$$

The main result of this section is the following theorem.

\begin{theorem}\label{thm: lattice ergodic average}
Let $T$ be a contractively regular operator on $L^p(\Omega)$ with $1<p<\infty$. Let $2<q<\infty$ and $\B$ be an UMD lattice having Fatou property. Then  there exist a constant $C_{p,q}$ such that
\begin{align}\label{lattice ergodic average}
\|(M_n(T)f)_{n\geq0}\|_{L^p(\B(v_q))}\|\leq C_{p,q}\|\|f\|\|_{p},\;\forall f\in L^p(\Omega;\B).
\end{align}
\end{theorem}

This proof is based on the transference principle, together with Theorem \ref{thm: lattice differential operator on Z1} . 

\begin{proof}
By the dilation theorem for regular operator (see e.g. \cite{Pel76}), There exists another measure space $(\hat{\Omega},\hat{\mu})$, two positive contractions $J:L^p(\Omega)\rightarrow L^p(\hat{\Omega})$ and $Q:L^p(\hat{\Omega})\rightarrow L^p(\Omega)$ and an isometric invertible operator $U: L^p(\hat{\Omega})\rightarrow L^p(\hat{\Omega})$ such that
$$T^k=QU^kJ,\;k\geq0.$$
Moreover $U$ can be chosen so that $U$ and $U^{-1}$ are both contractively regular. Thus $\|U^{j}\|_r=1$ for any $j\in \mathbb{Z}$.

Given $f\in L^p(\Omega;\B)$. Since $\B$ has Fatou property, it is easy to check that $$\|(M_n(T)f)_{n\geq0}\|_{L^p(\B(v_q))}=\lim_{N\rightarrow\infty }\|(M_n(T)f)_{0\leq n\leq N}\|_{L^p(\B(v_q))}.$$
Now fix an integer $N\geq1$. For any $n\geq0$, we clearly have
$$M_n(T)=QM_n(U)J.$$
Since $\|Q\|_r\leq1$, it follows from (\ref{extension of T}) that
\begin{align}\label{intermediate estimate}
\|(M_n(T)f)_{0\leq n\leq N}\|_{L^p(\Omega;\B(v_q))}\leq\|(M_n(U)J(f))_{0\leq n\leq N}\|_{L^p(\hat{\Omega};\B(v_q))}.
\end{align}
Using the regularity of $U^{-1}$, we have 
\begin{align*}
\|(M_n(U)J(f))_{0\leq n\leq N}\|_{L^p(\hat{\Omega};\B(v_q))}&=\|(M_n(U)U^{-\ell}U^{\ell}J(f))_{0\leq n\leq N}\|_{L^p(\hat{\Omega};\B(v_q))}\\
&\leq\|(M_n(U)U^{\ell}J(f))_{0\leq n\leq N}\|_{L^p(\hat{\Omega};\B(v_q))}
\end{align*}
for any integer $\ell$. 
Hence for any positive integer $L$, we have
\begin{align*}
&\|(M_n(U)J(f))_{0\leq n\leq N}\|^p_{L^p(\hat{\Omega};\B(v_q))}\\
&\leq \frac{1}{L}\sum^L_{\ell=0}\|(M_n(U)U^{\ell}J(f))_{0\leq n\leq N}\|^p_{L^p(\hat{\Omega};\B(v_q))}\\
&=\frac{1}{L}\int_{\hat{\Omega}}\sum^L_{\ell=0}\|(M_n(U)U^{\ell}J(f))_{0\leq n\leq N}\|^p_{\B(v_q))}d\hat{\omega}.
\end{align*}
Now we define a $\B$-valued function on $\mathbb{Z}\times \hat{\Omega}$, $g(k,\hat{\omega})=\chi_{0, N+L}(k)U^{k}J(f)(\hat{\omega})$. Then
\begin{align*}
M_n(U)U^{\ell}J(f)(\hat{\omega})&=\frac{1}{n+1}\sum^{n}_{k=0}U^{k+\ell}J(f)(\hat{\omega})=A_n(g(\cdot,\hat{\omega}))(\ell).
\end{align*}
Now apply Theorem \ref{thm: lattice differential operator on Z1},  using the regularity of $U$ and $J$, we have
\begin{align*}
&\|(M_n(U)J(f))_{0\leq n\leq N}\|^p_{L^p(\hat{\Omega};\B(v_q))}\leq C^p_{p,q}\frac{1}{L}\int_{\hat{\Omega}}\sum^{N+L}_{\ell=0}\|g(\ell,\hat{\omega})\|_{\B}^p d\hat{\omega}\\
&\leq \leq C^p_{p,q}\frac{1}{L}\int_{\hat{\Omega}}\sum^{N+L}_{\ell=0}\|J(f)(\hat{\omega})\|_{\B}^p d\hat{\omega}\leq\frac{N+L}{L}C^p_{p,q}\|f\|^p_{L^p(\B)}.
\end{align*}
\end{proof}

A continuous version of Theorem \ref{thm: lattice ergodic average} holds still true. Given a strongly continuous semigroup $T=(T_t)_{t\geq0}$ on $L^p(\Omega)$, we let
$$M_t(T)=\frac{1}{t}\int^t_0T_sds,\;t>0,$$
defined in the strong sense.

\begin{corollary}\label{cor: lattice ergodic average}
Let $T=(T_t)_{t\geq0}$ be a strongly continuous semigroup on $L^p(\Omega)$ and assume that $T_t: L^p(\Omega)\rightarrow L^p(\Omega)$ is contractively regular for any $t\geq0$. Let $2<q<\infty$, $\B$ be an UMD lattice having Fatou property  and let $f\in L^p(\Omega;\B)$. Then for a.e. $(\omega,\lambda)\in(\Omega\times\Sigma)$, the family $((M_t(T)f)(\omega,\lambda))_{t\geq0}$ belongs to $V_q$ and
$$\left\|(\omega,\lambda)\rightarrow\|((M_t(T)f)(\omega,\lambda))_{t\geq0}\|_{V_q}\right\|_p\leq C_{p,q}\|\|f\|\|_p.$$
\end{corollary}

This corollary can be proved in a similar way as Corollary \ref{cor: lattice differential operator on R1}.
\begin{remark}
Obviously, Corollary \ref{cor: lattice ergodic average} implies
$$\left\|(\omega,\lambda)\rightarrow\|((M_t(T)f)(\omega,\lambda))_{t\geq0}\|_{\infty}\right\|_p\leq C_{p}\|\|f\|\|_p$$
for any $f\in L^p(\Omega;\B)$, which is exactly Xu' Theorem 1 in \cite{Xu}.
\end{remark}

\section{Analytic semigroups}
In this section, we show that the inequality in Corollary \ref{cor: lattice ergodic average} holds still true when the ergodic averages are replaced by the semigroups themselves with an analyticity assumption. An operator $T$  on $L^p(\Omega)$ is called analytic, if $\sup_{n\geq0}n\|T^n-T^{n-1}\|<\infty$. We refer the reader to \cite{LeXu1} for more facts about analyticity. The following result is about the discrete contractively regular operators.

\begin{theorem}\label{thm: lattice analytic operator}
Let $1<p<\infty$,  $2<q<\infty$ and $\B$ be an UMD lattice having Fatou property.  Let $T$ be a contractively regular analytic operator on $L^p(\Omega)$. Then for any $f\in L^p(\Omega;\B)$ the family $(T^nf)_{n>0}$ belongs to $v_q$ for a.e. $(\omega,\lambda)\in\Omega\times\Sigma$ and we have an estimate
\begin{align}\label{lattice analytic semigroup operator}
\left\|(\omega,\lambda)\rightarrow \|(T^nf(\omega,\lambda))_{n>0}\|_{v_q}\right\|_{L^p(\Omega;\B)}\leq C_{p,q} \|\|f\|\|_p,\;\forall f\in L^p(\Omega;\B).
\end{align}
\end{theorem}

We will follow the strategy used in the scalar valued case \cite{LeXu2}, where the authors improved  a lot the original Stein's idea \cite{Ste70} to deal with the much stronger $q$-variaiton. And it turns out this new improvement is strong enough to enable us to deal with the vector-valued case. As in \cite{LeXu2}, we need the following lattice-valued square functions.

\begin{proposition}\label{pro: lattice littlewood-paley operator}
Let $1<p<\infty$ and $B$ be an UMD lattice. Let $T$ be a contractively regular analytic operator on $L^p(\Omega)$. Then
\begin{align}\label{lattice littlewood-paley operator}
\left\|\left(\sum^{\infty}_{n=0}({n+1})|T^{n+1}f-T^nf|^{2}\right)^{\frac{1}{2}}\right\|_{L^p(\B)}\leq C_{p}\|\|f\|\|_p,
\end{align}
with $C_p$ only depending on $p$ and $\B$.
\end{proposition}

This result will be proved after we show Theorem \ref{thm: lattice analytic operator}. With the square function estimates, the proof of Theorem \ref{thm: lattice analytic operator} can be done in a similar way as that in the scalar valued case.  We explain concisely the sketch of the proof, and refer the readers to Theorem 4.4 in \cite{LeXu2} for the details we omit here.

\begin{proof}
We begin with the following identity
\begin{align*}
T^{2n+1}&=\frac{1}{n}\sum^{2n}_{j=n}(j+1)T^j(T-I)-\frac{n+1}{n}(T^{2n+1}-T^n)\\
&+\frac{2n+1}{n}M_{2n}(T)-\frac{n+1}{n}M_n(T)\\
&=A_n-\frac{n+1}{n}B_n+\frac{2n+1}{n}M_{2n}(T)-\frac{n+1}{n}M_n(T),
\end{align*}
where
$$A_n=\frac{1}{n}\sum^{2n}_{j=n}(j+1)T^j(T-I),\;B_n=T^{2n+1}-T^n.$$
To finish the proof, it suffices to prove
\begin{equation}\label{estimate of An}
\|(A_n(f))_{n\geq0}\|_{L^p(\B(v_q))}\leq C_{p,q}\|\|f\|\|_p
\end{equation}
and
\begin{align}\label{estimate of Bn}
\|(B_n(f))_{n\geq0}\|_{L^p(\B(v_q))}\leq C_{p,q}\|\|f\|\|_p,
\end{align}
for all $f\in L^p(\Omega;\B)$. Since then by Lemma 4.3 in \cite{LeXu2} and Theorem \ref{thm: lattice ergodic average}, we have
\begin{align*}
&\|(\frac{n+1}{n}B_n(f))_{n\geq0}\|_{L^p(\B(v_q))}\\
&\leq 3\|(\frac{n+1}{n})_{n\geq0}\|_{v_1}\|(B_n(f))_{n\geq0}\|_{L^p(\B(v_q))}\leq C_{p,q}\|\|f\|\|_p,
\end{align*}
as well as
\begin{align*}
&\|(\frac{2n+1}{n}M_{2n}(T)(f))_{n\geq0}\|_{L^p(\B(v_q))},\;\|(\frac{n+1}{n}M_n(T)(f))_{n\geq0}\|_{L^p(\B(v_q))} \\
&\leq C_{p,q}\|\|f\|\|_p,
\end{align*}
and  finally
\begin{align*}
&\|(T^n(f))_{n\geq0}\|_{L^p(\B(v_q))}\\
&=\|(T^{2n+1}(f)-B_n(f))_{n\geq0}\|_{L^p(\B(v_q))}\leq C_{p,q}\|\|f\|\|_p,
\end{align*}
by using twice triangle inequalities.

Furthermore, using the same arguments in \cite{LeXu2}, for any fixed $(\omega,\sigma)\in \Omega\times\Sigma$, we can show
\begin{align*}
&\|(A_n(f)(\omega,\sigma))_{n\geq0}\|_{v_q},\;(B_n(f)(\omega,\sigma))_{n\geq0}\|_{v_q}\\
&\leq C\left(\sum^{\infty}_{n=0}({n+1})|(T^{n+1}f-T^nf)(\omega,\sigma)|^{2}\right)^{\frac{1}{2}}.
\end{align*}
Thus the desired estimates (\ref{estimate of An}) and (\ref{estimate of Bn}) follow from Proposition \ref{pro: lattice littlewood-paley operator}.
\end{proof}

A scalar valued versin of Proposition \ref{pro: lattice littlewood-paley operator} has been obtained in \cite{LeXu2}, where the authors used the theory of  $H^\infty$ functional calculus. The proof is very complicated. Actually, the proof of square function estimates constitutes the large portion of the paper \cite{LeXu1}. Recently, Xu \cite{Xu} has established $H^\infty$ functional calculus for analytic semigroups on vector-valued $L^p$ spaces, which allows us to adapt the idea in \cite{LeXu2} to the vector-valued square function estimates (\ref{pro: lattice littlewood-paley operator}). Let us give a sketch of the proof. We refer the readers to Theorem 3.3 in \cite{LeXu1} for the details we omit here.

\begin{proof}
Define a sequence of functions on $\mathbb{C}$, $F_n(z)=n^{1/2}z^{n-1}(z-1)$. Its $\gamma$-norm is defined as
$$\|(F_n)_{n\geq1}\|_{\gamma}=\sup\{\Big(\sum_{n\geq1}|F_n(z)|^2\Big)^{1/2},\;z\in B_{\gamma}\}$$
where $B_{\gamma}$ is the convex hull of 1 and the disc $D(0,\sin\gamma)$ with $\gamma$ being any real number $w(A)<\gamma<\pi/2$ and $w(A)$ being the type of $A=I-T$.
It is easy to check $\|(F_n)_{n\geq1}\|_{\gamma}<\infty$ (see e.g. \cite{LeXu2}).
The aim is to show
\begin{align*}
\|\Big(\sum_{n\geq1}|F_n(T)f|^2\Big)^{1/2}\|_{L^p(\B)}\leq C_p\|(F_n)_n\|_\gamma\|f\|_{L^p(\B)}.
\end{align*}
Define $G_n(z)=F_n(1-z)$. Since $1<p<\infty$ and $\B$ being UMD whence reflexive, by duality
\begin{align*}
\|\Big(\sum_{n\geq1}|F_n(T)f|^2\Big)^{1/2}\|_{L^p(\B)}=\sup\{|\sum_{n\geq1}\langle G_n(A)f,g_n\rangle|,\;(g_n)_n\in L^{p'}(\B^*(\ell^2))\}.
\end{align*}
Using Dunford functional calculus, we have
$$G_n(A)=\frac{1}{\pi i}\int_{L_\gamma}G_n(\lambda)A(\lambda-A)^{-1}(\lambda+A)^{-1}d\lambda$$
where $L_{\gamma}$ is the boundary of $1-B_\gamma$ oriented counterclockwise.
Hence
\begin{align*}
\langle G_n(A)f,g_n\rangle&=\frac{1}{\pi i}\int_{L_{\gamma}}G_n(\lambda)\langle A(\lambda-A)^{-1}(\lambda+A)^{-1}f,g_n\rangle d\lambda\\
&=\frac{1}{\pi i}\int_{L_{\gamma}}G_n(\lambda)\langle A^{1/2}(\lambda-A)^{-1}f,(A^*)^{1/2}(\lambda+A^*)^{-1}g_n\rangle d\lambda.
\end{align*}
Then by H\"older inequality
\begin{align*}
|\sum_{n\geq1}\langle G_n(A)f,g_n\rangle|&\leq \frac1\pi\|\Big(\int_{L_\gamma}\sum_n|G_n(\lambda)A^{1/2}(\lambda-A)^{-1}f|^2|d\lambda|\Big)^{1/2}\|_{L^p(\B)}\\
&\times \|\Big(\int_{L_\gamma}\sum_n|(A^*)^{1/2}(\lambda+A^*)^{-1}g_n|^2|d\lambda|\Big)^{1/2}\|_{L^{p'}(\B^*)}\\
&=\frac{1}{\pi}I\times II.
\end{align*}
The term $I$ can be further estimated as 
\begin{align*}
I&\leq\sup_{\lambda\in L_\gamma}(\sum_{n\geq1}|G_n(\lambda)|^2)^{1/2}\|\Big(\int_{L_\gamma}|A^{1/2}(\lambda-A)^{-1}f|^2|d\lambda|\Big)^{1/2}\|_{L^p(\B)}\\
&=\|(F_n)_n\|_{\gamma}\|\Big(\int_{L_\gamma}|A^{1/2}(\lambda-A)^{-1}f|^2|d\lambda|\Big)^{1/2}\|_{L^p(\B)}.
\end{align*}
Thus we are reduced to prove
\begin{align}\label{estimate of pb}
\|\Big(\int_{L_\gamma}|A^{1/2}(\lambda-A)^{-1}f|^2|d\lambda|\Big)^{1/2}\|_{L^p(\B)}\leq C_{p,q}\|f\|_{L^p(\B)}
\end{align}
and $||\leq C_{p',q}\|(g_n)_n\|_{L^{p'}(\B^*(\ell^2))}$, which can be also reduced to the estimates of the form (\ref{estimate of pb}) using Khintchine's inequality and Fubini's theorem (see \cite{LeXu1})
\begin{align}\label{estimate of pbstar}
\|\Big(\int_{L_\gamma}|(A^*)^{1/2}(\lambda+A^*)^{-1}g|^2|d\lambda|\Big)^{1/2}\|_{L^{p'}(\B^*)}\leq C_{p,q}\|g\|_{L^{p'}(\B)}.
\end{align}
 
The two estimates (\ref{estimate of pb}) and (\ref{estimate of pbstar}) are the vector-valued version of (3.4) and (3.5) in \cite{LeXu1}, which can be obtained by the same arguments in \cite{LeXu1} but using vector-valued results of $H^\infty$ functional calculus established recently in \cite{Xu}.

\end{proof}

\begin{remark}
If we define the $m$-difference sequence of $(T^n)_{n\geq0}$ for $T: L^p(\Omega)\rightarrow L^p(\Omega)$, $(\Delta^m_n)_{n\geq0}$ as
$$\Delta^m_n\equiv\Delta^m_n(T)=T^n(T-I)^m,$$
we then have a $m$ order version of square functions estimates using similar arguments (see \cite{LeXu1} and \cite{Xu})
\begin{align*}
\left\|\left(\sum^{\infty}_{n=0}\frac{1}{n+1}|(n+1)^{m+1}\Delta^{m+1}_n(f)|^{2}\right)^{\frac{1}{2}}\right\|_{L^p(\B)}\leq C_{p,m}\|\|f\|\|_p.
\end{align*}
Thus using similar arguments in the proof of Theorem \ref{thm: lattice analytic operator} (see \cite{LeXu2}), we obtain $m$ order version of Theorem \ref{thm: lattice analytic operator}
\begin{align}\label{m order variation}
\|(n^m\Delta^m_n(f))_{n\geq1}\|_{L^p(\B(v_q))}\leq C_{p,q,m} \|\|f\|\|_p,\;\forall f\in L^p(\Omega;\B).
\end{align}
\end{remark}

Now let us state the variational inequality for the continuous analytic semigroup.
\begin{corollary}\label{cor: lattice analytic semigroup}
Let $1<p<\infty$, $2<q<\infty$ and $\B$ be an UMD lattice having Fatou property. Let $(T_t)_{t\geq0}$ be a bounded analytic semigroup on $L^p(\Omega)$ and assume that $T_t$ is a contractively regular for any $t\geq0$. Then for any $f\in L^p(\Omega;\B)$, the family $(T_tf)_{t>0}$ belongs to $V_q$ for a.e. $(\omega,\lambda)\in\Omega\times\Sigma$ and we have an estimate
\begin{align}\label{lattice analytic semigroup}
\left\|(\omega,\lambda)\rightarrow \|(T_tf(\omega,\lambda))_{t>0}\|_{V_q}\right\|_{L^p(\Omega;\B)}\leq C_{p,q}\|\|f\|\|_p,\;\forall f\in L^p(\Omega;\B).
\end{align}
\end{corollary}

\begin{remark}
A vector-valued version of the sectorial maximal inequality has recently been proved by Xu \cite{Xu}.
\end{remark}

\begin{proof}
We will prove the variational inequality of order $m$ as in the discrete case (\ref{m order variation}).
Let $m\geq0$ be an integer. It follows from the lemma in Page 72 of \cite{Ste70} that the function
$$t\rightarrow t^m\frac{\partial^m}{\partial t^m}(T_t(f))(\omega,\lambda)$$
is continuous  for a.e. $(\omega,\lambda)\in\Omega\times\Sigma$.  To finish the proof, it suffices to show that  $T_t$'s satisfy the variational estimates (\ref{m order variation}) uniformly in $t$. Indeed, then using an approximation argument as in the proof of Corollary 4.2 in \cite{LeXu2}, we deduce that for any $0<t_0<t_1<\dotsm<t_N$, we have
\begin{align*}
&\|t^m\frac{\partial^m}{\partial t^m}(T_t(f))|_{t_0}, t^m\frac{\partial^m}{\partial t^m}(T_t(f))|_{t_1},\dotsm, t^m\frac{\partial^m}{\partial t^m}(T_t(f))|_{t_N}\|_{L^p(\B(v_q))}\\
&\leq C_{p,q}\|\|f\|\|_p.
\end{align*}
The desired result then follows from Lemma \ref{lem: approximation properties 4}.

Note that $T_t$'s are analytic uniformly, hence satisfy the estimate (\ref{lattice littlewood-paley operator}) uniformly. And the proof is finished. However, in the following we prefer to give another proof, which does not depend on Proposition \ref{pro: lattice littlewood-paley operator}. We shall deduce the estimate (\ref{lattice littlewood-paley operator})  for $T_t$'s directly from Xu's  square function inequality (see Proposition 10 of \cite{Xu})
\begin{align}\label{lattice littlewood-paley continuous}
\left\|\big(\int^{\infty}_0|s^{m+1}\frac{\partial^{m+1}T_s(f)}{\partial t^{m+1}}|^2\frac{ds}{s}\big)^{1/2}\right\|_{L^p(\B)}\lesssim\|\|f\|\|_p.
\end{align}
It is each to check that the estimate (\ref{lattice littlewood-paley operator}) is equivalent to the following inequality modulo a constant depending on $m$,
\begin{equation}\label{lattice littlewood-paley t}
\left\|\left(\sum^{\infty}_{n=0}\frac{1}n|n^{m+1}\Delta^{m+1}_n(T_t)(f)|^{2}\right)^{\frac{1}{2}}\right\|_{L^p(\B)}\leq C_{m,p}\|\|f\|\|_p.
\end{equation}
Estimates (\ref{lattice littlewood-paley t}) will be dealt with in two cases $0<t<\varepsilon$ and $t\geq\varepsilon$ for some sufficiently small $\varepsilon$. We first deal with the case $0<t<\varepsilon$. Since $(T_t)_{t\geq0}$ is an analytic semigroup of regular contractions, and $B$ an UMD lattice, Theorem 4 in \cite{Xu} yields the generator $-A$ of the semigroup admits a bounded $H^{\infty}(\Sigma_{\sigma})$ functional calculus for some $\sigma<\pi/2$. Hence taking $\varphi(z)=z^{m+1}e^{-z}$ which is a bounded analytic function,  then
$$\varphi(sA)=(sA)^{m+1}e^{-sA}=(-1)^{m+1}s^{m+1}\frac{\partial^{m+1}T_s(f)}{\partial s^{m+1}}.$$
Then by the lemma in page 72 of \cite{Ste70} and its proof, we can find very small $\varepsilon$ such that for any $0<t<\varepsilon$, the left hand side of inequality (\ref{lattice littlewood-paley continuous}) is equivalent to the quantity
\begin{align*}
\left\|\big(\sum^{\infty}_{n=0}n^{2m+1}t^{2m+2}|\frac{\partial^{m+1}T_s(f)}{\partial s^{m+1}}|_{s=nt}|^2\big)^{1/2}\right\|_{L^p(\B)}
\end{align*}
which is further equivalent to the left hand side of inequality (\ref{lattice littlewood-paley t})  for small $t$ since the partial derivative insider the sum can be approximated by the difference quotient
$$\frac{\Delta^{m+1}_n(T_t)}{t^{m+1}}.$$
Estimates (\ref{lattice littlewood-paley t}) in the case $t\geq\varepsilon$ follows from the case $0<t<\varepsilon$ by the following observation
\begin{align}\label{observation}
\left(\sum^{\infty}_{n=0}n^{2m+1}|\Delta^{m+1}_n(T_{2t})(f)|^{2}\right)^{\frac{1}{2}}\leq\left(\sum^{\infty}_{n=0}n^{2m+1}|\Delta^{m+1}_n(T_{t})(f)|^{2}\right)^{\frac{1}{2}}.
\end{align}
Let us prove this observation, first of all, we have the following identity
\begin{align}\label{identity}
\Delta^{m+1}_n(T_{2t})=\sum^{m+1}_{k=0}\begin{pmatrix} m+1\\ k\end{pmatrix}\Delta^{m+1}_{2n+k}(T_t),
\end{align}
which can be easily checked as follows by definition
\begin{align*}
\sum^{m+1}_{k=0}\Delta^{m+1}_{2n+k}(T_t)&=\sum^{m+1}_{k=0}\begin{pmatrix} m+1\\ k\end{pmatrix}\sum^{m+1}_{j=0} \begin{pmatrix} m+1\\ j\end{pmatrix}(-1)^{m+1-j}T_t^{2n+k+j}\\
&=T_t^{2n}(T+I)^{m+1}(T-I)^{m+1}=\Delta^{m+1}_n(T_{2t}).
\end{align*}
By triangle inequality and Jensen inequality, the identity (\ref{identity}) implies the left hand side of inequality (\ref{observation}) is smaller than
\begin{align*}
&\sum^{m+1}_{k=0}\begin{pmatrix} m+1\\ k\end{pmatrix}\left(\sum^{\infty}_{n=0}n^{2m+1}|\Delta^{m+1}_{2n+k}(T_{t})(f)|^{2}\right)^{\frac{1}{2}}\\
&\leq\sum^{m+1}_{k=0}\begin{pmatrix} m+1\\ k\end{pmatrix}\left(\frac{1}{2^{2m+1}}\sum^{\infty}_{n=0}{(2n+k)}^{2m+1}|\Delta^{m+1}_{2n+k}(T_{t})(f)|^{2}\right)^{\frac{1}{2}}\\
&\leq\left(\frac{1}{2^{2m+1}}\sum^{m+1}_{k=0}\begin{pmatrix} m+1\\ k\end{pmatrix}\sum^{m+1}_{k=0}\begin{pmatrix} m+1\\ k\end{pmatrix}\sum^{\infty}_{n=0}{(2n+k)}^{2m+1}|\Delta^{m+1}_{2n+k}(T_{t})(f)|^{2}\right)^{\frac{1}{2}}\\
&\leq\left(\frac{1}{2^{2m+1}}2^{m+1}\begin{pmatrix} m+1\\ \big[\frac{m+1}{2}\big]\end{pmatrix}\sum^{m+1}_{k=0}\sum^{\infty}_{n=0}{(2n+k)}^{2m+1}|\Delta^{m+1}_{2n+k}(T_{t})(f)|^{2}\right)^{\frac{1}{2}},
\end{align*}
which is controlled by the right hand side of (\ref{observation}) due to the fact
$$\begin{pmatrix} m+1\\ \big[\frac{m+1}{2}\big]\end{pmatrix}\leq 2^{m},\;\;\forall m\geq0.$$
\end{proof}

\section{Applications}
In this section, we discuss about the consequences of the variational inequalities.
Let $T=(T_t)_{t\geq0}$ be a strongly continuous semigroup of contractively regular operators on $L^p(\Omega)$.  As argued by Taggart \cite{Tag09},  for any Banach space $\B$ and any $f\in L^p(\Omega;B)$, $T_tf\rightarrow f$ in the $L^p(B)$-norm as $t\rightarrow0^+$. This implies $M_t(T)f\rightarrow f$ in the $L^p(B)$-norm as $t\rightarrow0^+$.

Let $A$ denote the infinitesimal generator of $T$. By mean ergodic theorem, we have direct sum decomposition
$$L^p(\Omega)=N(A)\oplus \overline{R(A)}.$$
If we let $P_A: L^p(\Omega)\rightarrow L^p(\Omega)$ denotes the corresponding projection onto $N(A)$, then
$$M_t(T)f\rightarrow P_A(f),\;\;\mathrm{as}\;t\rightarrow\infty$$
and if further $(T_t)_{t\geq0}$ is a bounded analytic semigroup, then
$$T_tf\rightarrow P_A(f),\;\;\mathrm{as}\;t\rightarrow\infty$$
in $L^p$-norm for any $f\in L^p(\Omega)$. By similar reason, these limits hold also true for any $f\in L^p(\Omega;\B)$.

Now applying Corollary \ref{cor: lattice ergodic average} and Theorem \ref{thm: lattice analytic semigroup}, we deduce the following individual ergodic theorems.
\begin{corollary}
Let $T=(T_t)_{t\geq0}$ be a strongly continuous semigroup of contractively regular operators on $L^p(\Omega)$. Let  $f\in L^p(\Omega;\B)$ with $\B$  an UMD lattice on $(\Sigma,\nu)$ having Fatou property, then for almost every $(\omega,\lambda)\in\Omega\times\Sigma$,
$$[M_t(T)f](\omega,\lambda)\rightarrow [P_A(f)](\omega,\lambda) ,\;\;\mathrm{as}\;t\rightarrow\infty$$
and
$$[M_t(T)f](\omega,\lambda)\rightarrow f(\omega,\lambda) ,\;\;\mathrm{as}\;t\rightarrow0^+.$$
\end{corollary}
\begin{corollary}
Let $T=(T_t)_{t\geq0}$ be a bounded analytic semigroup on $L^p(\Omega)$ and assume that $T_t$ is contractively regular for any $t\geq0$. Let  $f\in L^p(\Omega;\B)$ with $\B$ an UMD lattice on $(\Sigma,\nu)$ having Fatou property, then for almost every $(\omega,\lambda)\in\Omega\times\Sigma$,
$$[T_tf](\omega,\lambda)\rightarrow [P_A(f)](\omega,\lambda) ,\;\;\mathrm{as}\;t\rightarrow\infty$$
and
$$[T_tf](\omega,\lambda)\rightarrow f(\omega,\lambda) ,\;\;\mathrm{as}\;t\rightarrow0^+.$$
\end{corollary}

These pointwise convergences can be deduced from the corresponding maximal inequality established recently by Xu in \cite{Xu} via the Banach principle.  The point of this paper is that Corollary \ref{cor: lattice ergodic average} and Theorem \ref{thm: lattice analytic semigroup} provides us quantitative information on the rate of the convergence.  This requires the notion of $\lambda$-jump functions. For any $\lambda>0$ and any family $a=(a_t)_{t\geq0}$ of complex numbers. One defines $N(a,\lambda)$ to be the supremum of all integers $N\geq0$ for which there is an increasing sequence
$$0<s_1<t_1\leq s_2<t_2\leq\dotsm\leq s_N<t_N$$
so that $|a_{t_k}-a_{s_k}|>\lambda$ for each $k=1,\dotsm,N$. It is clear that for any $1\leq q<\infty$,
$$\lambda^qN(a,\lambda)\leq \|a\|^q_{V_q}.$$

By Corollary \ref{cor: lattice ergodic average} and Theorem \ref{thm: lattice analytic semigroup}, we immediately obtain the following jump estimates.
\begin{corollary}
Let $T=(T_t)_{t\geq0}$ be a strongly continuous semigroup of contractively regular operators on $L^p(\Omega)$. Let  $f\in L^p(\Omega;\B)$ with $\B$ an UMD lattice on $(\Sigma,\nu)$ having Fatou property and $2<q<\infty$, then we have
$$\left\|\omega\rightarrow N\big(([M_t(T)f](\omega))_{t\geq0},\lambda\big)^{\frac{1}{q}}\right\|_{L^p(\B)}\lesssim\frac{\|f\|_{L^p(\B)}}{\lambda},$$
and for any $K>0$, we also have
$$\mu\times\nu\big\{\omega\in\Omega|N\big(([M_t(T)f](\omega))_{t\geq0},\lambda\big)>K\big\}\lesssim\frac{\|f\|_{L^p(B)}}{\lambda^pK^{\frac{p}{q}}} .$$
\end{corollary}

\begin{corollary}
Let $T=(T_t)_{t\geq0}$ be a bounded analytic semigroup on $L^p(\Omega)$ and assume that $T_t$ is contractively regular for any $t\geq0$. Let  $f\in L^p(\Omega;\B)$ with $\B$ be an UMD lattice on $(\Sigma,\nu)$ having Fatou property and $2<q<\infty$, then we have
$$\left\|\omega\rightarrow N\big(([T_tf](\omega))_{t\geq0},\lambda\big)^{\frac{1}{q}}\right\|_{L^p(B)}\lesssim\frac{\|f\|_{L^p(\B)}}{\lambda},$$
and for any $K>0$, we also have
$$\mu\times\nu\big\{\omega\in\Omega|N\big(([T_tf](\omega))_{t\geq0},\lambda\big)>K\big\}\lesssim\frac{\|f\|_{L^p(\B)}}{\lambda^pK^{\frac{p}{q}}} .$$
\end{corollary}

As in the scalar-valued case (see e.g. \cite{Bou89}, \cite{JSW08}), we may expect that the jump estimate would be true when $q=2$. This fact will be proved and appear elsewhere.

\vskip3pt

\noindent \textbf{Acknowledgement.} Guixiang Hong is supported by MINECO: ICMAT Severo Ochoa project SEV-2011-0087 and ERC Grant StG-256997-CZOSQP (EU). Tao Ma is partially supported by NSFC No. 11271292.

\vskip30pt

\end{document}